\documentclass{elsart}
\usepackage{amssymb}
\usepackage{graphicx}

\begin{document}

\begin{frontmatter}

\title{The~general Li\'{e}nard~polynomial~system\thanksref{label1}}
\thanks[label1]{This work was supported by the Netherlands Organization for Scientific Research (NWO).
The author is also very grateful to the Max Planck Institute for Ma\-the\-matics (Bonn) and the Johann
Bernoulli Institute for Mathematics and Computer Science (Groningen) for hospitality and support during
his stay in 2011--\,2012.
}

\author{Valery A. Gaiko}

\ead{valery.gaiko@yahoo.com}

\address{National~Academy~of~Sciences~of~Belarus,~United~Institute~of~Informatics~Problems,
Surganov~Str.~6,~Minsk~220012,~Belarus}

\begin{abstract}
In this paper, applying a canonical system with field rotation parameters and using geometric
properties of the spirals filling the interior and exterior domains of limit cycles, we solve
first the problem on the maximum number of limit cycles surrounding a unique singular point for
an arbitrary polynomial system. Then, by~means of the same bifurcationally geometric approach,
we solve the limit cycle problem for a general Li\'{e}nard polynomial system with an arbitrary
(but finite) number of singular points. This is related to the solution of Hil\-bert's
sixteenth problem on the maximum number and relative position of limit cycles for
planar polynomial dynamical systems.
\end{abstract}

\end{frontmatter}

\section{Introduction}

We will consider Li\'{e}nard equations of the form
    $$
    \ddot{x}+f(x)\,\dot{x}+g(x)=0.
    \eqno(1.1)
    $$
In the phase plane, the representation of equation (1.1) is given by the dynamical system
    $$
    \dot{x}=y, \quad \dot{y}=-g(x)-f(x)y.
    \eqno(1.2)
    $$
    \par
There are many examples in the natural sciences and technology in which this and related systems
are applied \cite{BL}--\cite{sml}. Such systems are often used to model either mechanical or electrical,
or biomedical systems, and in the literature, many systems are transformed into Li\'{e}nard type to aid
in the investigations. They can be used, e.\,g., in certain mechanical systems, where $f(x)$ represents
a coefficient of the damping force and $g(x)$ represents the restoring force or stiffness, when modeling
wind rock phenomena and surge in jet engines~\cite{agan},~\cite{ochbc}. Such systems can be also used
to model resistor-inductor-capacitor circuits with non-linear circuit elements. Recently, e.\,g., the
Li\'{e}nard system (1.2) has been shown to describe the operation of an optoelectronics circuit that
uses a resonant tunnelling diode to drive a laser diode to make an optoelectronic vol\-tage controlled
oscillator \cite{srlfwi}. There are also some examples of using Li\'{e}nard type systems in ecology
and epidemiology \cite{mrr}.
    \par
In this paper, we suppose that system (1.2), where $f(x)$ and $g(x)$ are arbitrary polynomials of~$x,$
has an anti-saddle (a node or a focus, or a center) at the origin and write it in the form
    $$
    \dot{x}=y,
    \quad
    \dot{y}=-x\,(1+\beta_{1}\,x+\ldots+\beta_{2l}\,x^{2l})
    +y\,(\alpha_{0}+\alpha_{1}\,x+\ldots+\alpha_{2k}\,x^{2k}).
    \eqno(1.3)
    $$
    \par
Note that for $g(x)\equiv~x,$ by the change of variables $X=x$ and $Y=y+F(x),$ where
$F(x)=\int_{0}^{x}\!f(s)\,\textrm{d}s,$ (1.3) is reduced to an equivalent system
    $$
    \dot{X}=Y-F(X), \quad \dot{Y}=-X
    \eqno(1.4)
    $$
which can be written in the form
    $$
    \dot{x}=y, \quad \dot{y}=-x+F(y)\\[-4mm]
    \eqno(1.5)
    $$
or\\[-2mm]
    $$
    \dot{x}=y,
    \quad
    \dot{y}=-x+\gamma_{1}\,y+\gamma_{2}\,y^{2}+\gamma_{3}\,y^{3}+\ldots
    +\gamma_{2k}\,y^{2k}+\gamma_{2k+1}\,y^{2k+1}.\\[1mm]
    \eqno(1.6)
    $$
    \par
Therefore, we can conclude that our previous results \cite{g08a}, \cite{g09a} agree with the conjecture
of~\cite{ldmp} on the maximum number of limit cycles for the classical Li\'{e}nard polynomial system~(1.6).
In~\cite{g08a},~\cite{g09a}, we have presented a solution of Smale's thirteenth problem~\cite{sml} proving
that the Li\'{e}nard system~(1.6) with a polynomial of degree $2k+1$ can have at most $k$ limit cycles.
In \cite{g97}--\cite{g09b}, we have also presented a solution of Hilbert's sixteenth problem in the quadratic
case of polynomial systems proving that for quadratic systems four is really the maximum number of limit cycles
and $(3\!:\!1)$ is their only possible distribution. We have established some preliminary results on generalizing
our ideas and methods to special cubic, quartic and other polynomial dynamical systems as~well.
In \cite{gvh04}, e.\,g., we have constructed a canonical cubic dynamical system of Kukles type and
have carried out the global qualitative analysis of its special case corresponding to a generalized
Li\'{e}nard equation. In particular, it has been shown that the foci of such a Li\'{e}nard system can be
at most of second order and that such system can have at most three limit cycles in the whole phase plane.
Moreover, unlike all previous works on the Kukles-type systems, global bifurcations of limit and separatrix
cycles using arbitrary (including as large as possible) field rotation parameters of the canonical system
have been studied. As a result, a classification of all possible types of separatrix cycles for the
generalized Li\'{e}nard system has been obtained and all possible distributions of its limit cycles have
been found. In~\cite{gvh09a},~\cite{gvh09b}, we have completed the global qualitative analysis of a planar
Li\'{e}nard-type dynamical system with a piecewise linear function containing an arbitrary number of dropping
sections and approximating an arbitrary polynomial function. In~\cite{bnrs1}--\cite{brg}, we have also completed
the global qualitative analysis of a quartic dynamical system which models the dynamics of the populations
of predators and their prey in a given ecological system. In~\cite{g11}, we have studied of multiple limit
cycle bifurcations in the well-known FitzHgh--Nagumo neuronal model.
    \par
We use the obtained results and develop our methods for studying limit cycle bifurcations of polynomial dynamical
systems in this paper as well. In~Section~2, applying a canonical system with field rotation parameters and using
geometric properties of the spirals filling the interior and exterior domains of limit cycles, we solve first
the problem on the maximum number of limit cycles surrounding a unique singular point for an arbitrary
polynomial system. Then, in Section~3, by means of the same bifurcationally geometric approach, we solve
the limit cycle problem for the general Li\'{e}nard polynomial system~(1.3) with an arbitrary (but finite)
number of singular points proving that system (1.3) can have at most $k+l$~limit cycles, $k$~surrounding
the origin and $l$~surrounding one by one the other singularities of~(1.3). This is related to the solution
of Hil\-bert's sixteenth problem on the maximum number and relative position of limit cycles for planar
polynomial dynamical systems.

\section{Limit cycles surrounding a unique singular point}

Consider first an arbitrary dynamical system
    $$
    \dot{x}=P_{n}(x,y,\mu_{1},\ldots,\mu_{k}),
    \quad
    \dot{y}=Q_{n}(x,y,\mu_{1},\ldots,\mu_{k}),
    \eqno(2.1)
    $$
where $P_{n}$ and $Q_{n}$ are polynomials in the real variables
$x,\:y$ and not greater than $n$-th degree containing $k$ field
rotation parameters, $\mu_{1},\ldots,\mu_{k},$ and having an
anti-saddle at the origin. We prove the following theorem.
    \par
    \medskip
\noindent\textbf{Theorem 2.1.}
    \emph{The polynomial system $(2.1)$ containing $k$ field rotation parameters
and having a singular point of center type at the origin for the zero values of these
parameters can have at most $k-1$~limit cycles surrounding the origin.}
    \par
\noindent\textbf{Proof.} Let all the parameters of (2.1) vanish and suppose that
the obtained system
    $$
    \dot{x}=P_{n}(x,y,0,\ldots,0),
    \quad
    \dot{y}=Q_{n}(x,y,0,\ldots,0)\\[2mm]
    \eqno(2.2)
    $$
has a singular point of center type at the origin.
    \par
Input successively the field rotation parameters, $\mu_{1},\ldots,\mu_{k},$
into this system (see \cite{g08a}, \cite{g09a}). Suppose, e.\,g., that
$\mu_{1}>0$ and that the vector field of the system
    $$
    \dot{x}=P_{n}(x,y,\mu_{1},0,\ldots,0),
    \quad
    \dot{y}=Q_{n}(x,y,\mu_{1},0,\ldots,0)\\[2mm]
    \eqno(2.3)
    $$
is rotated counterclockwise turning the origin into a stable focus under
increasing the parameter $\mu_{1}$ \cite{BL}, \cite{g08a}, \cite{g09a}.
    \par
Fix $\mu_{1}$ and input the parameter $\mu_{2}$ into (2.3)
changing it so that the field of the system
    $$
    \dot{x}=P_{n}(x,y,\mu_{1},\mu_{2},0,\ldots,0),
    \quad
    \dot{y}=Q_{n}(x,y,\mu_{1},\mu_{2},0,\ldots,0)
    \eqno(2.4)
    $$
would be rotated in the opposite direction (clockwise). Suppose this occurs
for $\mu_{2}<0.$ Then, for some value of this parameter, a limit cycle
will appear in system~(2.4). There are three logical possibilities
for such a bifurcation: 1)~the limit cycle appears from the focus
at the origin; 2)~it can also appear from some separatrix cycle
surrounding the origin; 3)~the limit cycle appears from a
so-called ``trajectory concentration''. In the last case, the
limit cycle is semi-stable and, under further decreasing
$\mu_{2},$ it splits into two limit cycles (stable and unstable)
one of which then disappears at (or tends to) the origin and the
other disappears on (or tends to) some separatrix cycle
surrounding this point. But since the stability character of both
a~singular point and a~separatrix cycle is quite easily controlled~\cite{G},
this logical possibility can be excluded. Let us choose one of the two other
possibilities: e.\,g., the first one, the so-called Andronov--Hopf
bifurcation. Suppose that, for some value of~$\mu_{2},$ the focus at
the origin becomes non-rough, changes the character of its stability
and generates a stable limit cycle, $\Gamma\!_{1}.$
    \par
Under further decreasing $\mu_{2},$ three new logical
possibilities can arise: 1)~the limit cycle $\Gamma\!_{1}$
disappears on some separatrix cycle surrounding the origin;
2)~a~sepa\-ratrix cycle can be formed earlier than $\Gamma\!_{1}$
disappears on it and then it generates one more (unstable) limit
cycle, $\Gamma_{2},$ which joins with $\Gamma\!_{1}$ forming a
semi-stable limit cycle, $\Gamma\!_{12},$ disappearing in a
``trajectory concentration'' under further decreasing $\mu_{2};$
3)~in the domain $D_{1}$ outside the cycle $\Gamma\!_{1}$ or in
the domain $D_{2}$ inside $\Gamma\!_{1},$ a semi-stable limit
cycle appears from a ``trajectory concentration'' and then splits
into two limit cycles (logically, the appearance of such
semi-stable limit cycles can be repeated).
    \par
Let us consider the third case. It is clear that, under decreasing
$\mu_{2},$ a semi-stable limit cycle cannot appear in the domain
$D_{2},$ since the focus spirals filling this domain will untwist
and the distance between their coils will increase because of the
vector field rotation. By contradiction, we can prove that a
semi-stable limit cycle cannot appear in the domain $D_{1}.$
    \par
Suppose it appears in this domain for some values of the parameters
$\mu_{1}^{*}>0$ and $\mu_{2}^{*}<0.$ Return to system~(2.2) and
change the inputting order for the field rotation parameters.
Input first the parameter $\mu_{2}<0\!:$
    $$
    \dot{x}=P_{n}(x,y,\mu_{2},0,\ldots,0),
    \quad
    \dot{y}=Q_{n}(x,y,\mu_{2},0,\ldots,0).
    \eqno(2.5)
    $$
Fix it under $\mu_{2}=\mu_{2}^{*}.$ The vector field of~(2.5)
is rotated clockwise and the origin turns into an unstable focus.
Inputting the parameter $\mu_{1}>0$ into (2.5), we get again
system (2.4) the vector field of which is rotated counterclockwise.
Under this rotation, a stable limit cycle, $\Gamma\!_{1},$ will appear
from some separatrix cycle. The limit cycle $\Gamma\!_{1}$ will contract,
the outside spirals winding onto this cycle will untwist and the distance
between their coils will increase under increasing $\mu_{1}$ to the value
$\mu_{1}^{*}.$ It follows that there are no values of $\mu_{2}^{*}<0$ and
$\mu_{1}^{*}>0$ for which a semi-stable limit cycle could appear in the
domain~$D_{1}.$
    \par
The second logical possibility can be excluded by controlling the
stability character of the separatrix cycle~\cite{G}. Thus,
only the first possibility is valid, i.\,e., system~(2.4)
has at most one limit cycle.
    \par
Let system (2.4) have the unique limit cycle $\Gamma\!_{1}.$
Fix the parameters $\mu_{1}>0,$ $\mu_{2}<0$ and input the third
parameter, $\mu_{3}>0,$ into this system supposing that $\mu_{3}$
rotates its vector field counterclockwise:
    $$
    \dot{x}=P_{n}(x,y,\mu_{1},\mu_{2},\mu_{3},0,\ldots,0),
    \quad
    \dot{y}=Q_{n}(x,y,\mu_{1},\mu_{2},\mu_{3},0,\ldots,0).
    \eqno(2.6)
    $$
Here we can have two basic possibilities: 1)~the limit cycle
$\Gamma\!_{1}$ disappears at the origin; 2)~the second (unstable)
limit cycle, $\Gamma_{2},$ appears from the origin and, under
further increasing the parameter $\mu_{3},$ the cycle $\Gamma_{2}$
joins with $\Gamma\!_{1}$ forming a semi-stable limit cycle,
$\Gamma_{\!12},$ which disappears in a ``trajectory
concentration'' surrounding the origin. Besides, we can also
suggest that: 3)~in the domain $D_{2}$ bounded by the origin and
$\Gamma\!_{1},$ a semi-stable limit cycle, $\Gamma_{23},$ appears
from a ``trajectory concentration'', splits into an unstable
cycle, $\Gamma_{2},$ and a stable cycle, $\Gamma_{3},$ and then
the cycles $\Gamma\!_{1},$ $\Gamma_{2}$ disappear through a
semi-stable limit cycle, $\Gamma\!_{12},$ and the cycle
$\Gamma_{3}$ disappears through an Andronov--Hopf bifurcation;
4)~a~semi-stable limit cycle, $\Gamma_{34},$ appears in the domain
$D_{2}$ bounded by the cycles $\Gamma\!_{1},$ $\Gamma_{2}$ and,
for some set of values of the parameters, $\mu_{1}^{*},$
$\mu_{2}^{*},$ $\mu_{3}^{*},$ system (2.6) has at least
four limit cycles.
    \par
Let us consider the last, fourth case. It is clear that a
semi-stable limit cycle cannot appear either in the domain $D_{1}$
bounded on the inside by the cycle $\Gamma\!_{1}$ or in the domain
$D_{3}$ bounded by the origin and $\Gamma_{2}$ because of the
increasing distance between the spiral coils filling these domains
under increasing the parameter $\mu_{3}.$ To prove the impossibility
of the appearance of a semi-stable limit cycle in the domain
$D_{2},$ suppose the contrary, i.\,e., for some set of values of
the parameters, $\mu_{1}^{*}>0,$ $\mu_{2}^{*}<0,$ and
$\mu_{3}^{*}>0,$ such a semi-stable cycle exists. Return to
system~(2.2) again and input first the parameters $\mu_{3}>0,$
$\mu_{1}>0\!:$
    $$
    \dot{x}=P_{n}(x,y,\mu_{1},\mu_{3},0,\ldots,0),
    \quad
    \dot{y}=Q_{n}(x,y,\mu_{1},\mu_{3},0,\ldots,0).
    \eqno(2.7)
    $$
Fix these parameters under $\mu_{3}=\mu_{3}^{*},$
$\mu_{1}=\mu_{1}^{*}$ and input the parameter $\mu_{2}<0$ into~(2.7)
getting again system~(2.6). Since, by our assumption, this
system has two limit cycles for $\mu_{2}>\mu_{2}^{*},$ there
exists some value of the parameter, $\mu_{2}^{12}$
$(\mu_{2}^{*}<\mu_{2}^{12}<0),$ for which a semi-stable limit
cycle, $\Gamma\!_{12},$ appears in system (2.6) and then splits
into a stable cycle, $\Gamma\!_{1},$ and an unstable cycle,
$\Gamma_{2},$ under further decreasing $\mu_{2}.$ The formed
domain $D_{2}$ bounded by the limit cycles $\Gamma\!_{1},$
$\Gamma_{2}$ and filled by the spirals will enlarge, since, by the
properties of a field rotation parameter, the interior unstable
limit cycle $\Gamma_{2}$ will contract and the exterior stable
limit cycle $\Gamma\!_{1}$ will expand under decreasing $\mu_{2}.$
The distance between the spirals of the domain $D_{2}$ will
naturally increase, what will prevent the appearance of a
semi-stable limit cycle in this domain for $\mu_{2}<\mu_{2}^{12}.$
    \par
Thus, there are no such values of the parameters, $\mu_{1}^{*}>0,$
$\mu_{2}^{*}<0,$ $\mu_{3}^{*}>0,$ for which system (2.6)
would have an additional semi-stable limit cycle. Therefore, the
fourth case cannot be realized. The third case is considered absolutely
similarly. It follows from the first two cases that system~(2.6)
can have at most two limit cycles.
    \par
Suppose that system (2.6) has two limit cycles, $\Gamma\!_{1}$
and $\Gamma_{2},$ fix the parameters $\mu_{1}>0,$ $\mu_{2}<0,$
$\mu_{3}>0$ and input the fourth parameter, $\mu_{4}<0,$ into this
system supposing that $\mu_{4}$ rotates its vector field
clockwise:
    $$
    \dot{x}=P_{n}(x,y,\mu_{1},\ldots,\mu_{4},0,\ldots,0),
    \quad
    \dot{y}=Q_{n}(x,y,\mu_{1},\ldots,\mu_{4},0,\ldots,0).
    \eqno(2.8)
    $$
The most interesting logical possibility here is that when the
third (stable) limit cycle, $\Gamma_{3},$ appears from the origin
and then, under preservation of the cycles $\Gamma\!_{1}$ and
$\Gamma_{2},$ in the domain $D_{3}$ bounded on the inside by the
cycle $\Gamma_{3}$ and on the outside by the cycle $\Gamma_{2},$ a
semi-stable limit cycle, $\Gamma_{\!45},$ appears and then splits
into a stable cycle, $\Gamma\!_{4},$ and an unstable cycle,
$\Gamma_{5},$ i.\,e., when system~(2.8) for some set of values
of the parameters, $\mu_{1}^{*},$ $\mu_{2}^{*},$ $\mu_{3}^{*},$
$\mu_{4}^{*},$ has at least five limit cycles. Logically, such a
semi-stable limit cycle could also appear in the domain $D_{1}$
bounded on the inside by the cycle $\Gamma\!_{1},$ since, under
decreasing $\mu_{4},$ the spirals of the trajectories of (2.8)
will twist and the distance between their coils will decrease. On
the other hand, in the domain $D_{2}$ bounded on the inside by the
cycle $\Gamma_{2}$ and on the outside by the cycle $\Gamma\!_{1}$
and also in the domain $D_{4}$ bounded by the origin and
$\Gamma_{3},$ a semi-stable limit cycle cannot appear, since,
under decreasing $\mu_{4},$ the spirals will untwist and the
distance between their coils will increase. To prove the impossibility
of the appearance of a semi-stable limit cycle in the domains
$D_{3}$ and $D_{1},$ suppose the contrary, i.\,e., for some set
of values of the parameters, $\mu_{1}^{*}>0,$ $\mu_{2}^{*}<0,$
$\mu_{3}^{*}>0,$ and $\mu_{4}^{*}<0,$ such a semi-stable
cycle exists.
    \par
Return to system (2.2) again, input first the parameters
$\mu_{4}<0,$ $\mu_{2}<0$ and then the parameter $\mu_{1}>0\!:$
    $$
    \dot{x}=P_{n}(x,y,\mu_{1},\mu_{2},\mu_{4},0,\ldots,0),
    \quad
    \dot{y}=Q_{n}(x,y,\mu_{1},\mu_{2},\mu_{4},0,\ldots,0).
    \eqno(2.9)
    $$
Fix the parameters $\mu_{4},$ $\mu_{2}$ under the values
$\mu_{4}^{*},$ $\mu_{2}^{*},$ respectively. Under increasing
$\mu_{1},$ a separatrix cycle is formed around the origin
generating a stable limit cycle, $\Gamma\!_{1}.$ Fix $\mu_{1}$
under the value $\mu_{1}^{*}$ and input the parameter $\mu_{3}>0$
into (2.9) getting system~(2.8).
    \par
Since, by our assumption, system (2.8) has three limit cycles for
$\mu_{3}<\mu_{3}^{*},$ there exists some value of the parameter
$\mu_{3}^{23}$ $(0<\mu_{3}^{23}<\mu_{3}^{*})$ for which a
semi-stable limit cycle, $\Gamma_{23},$ appears in this system and
then splits into an unstable cycle, $\Gamma_{2},$ and a stable
cycle, $\Gamma_{3},$ under further increasing $\mu_{3}.$ The
formed domain $D_{3}$ bounded by the limit cycles $\Gamma_{2},$
$\Gamma_{3}$ and also the domain $D_{1}$ bounded on the inside
by the limit cycle $\Gamma\!_{1}$ will enlarge and the spirals
filling these domains will untwist excluding a possibility of
the appearance of a semi-stable limit cycle there.
    \par
All other combinations of the parameters $\mu_{1},$ $\mu_{2},$
$\mu_{3},$ and $\mu_{4}$ are considered in a similar way. It
follows that system (2.8) has at most three limit cycles. If we
continue the procedure of successive inputting the field rotation
parameters, $\mu_{5},$ $\mu_{6},\ldots,$ $\mu_{k},$ into system~(2.2),
it is possible to conclude that system (2.1) can have at most $k-1$~limit
cycles surrounding the origin. The theorem is proved.\qquad $\Box$

\section{Limit cycles of the general Li\'{e}nard polynomial system}

By means of the same bifurcationally geometric approach, we will consider
now the Li\'{e}nard polynomial system (1.3). The study of singular points
of system~(1.3) will use two index theorems by H.\,Poincar\'{e}, see \cite{BL}.
The definition of the Poincar\'{e} index is the following~\cite{BL}.
    \medskip
    \par
    \textbf{Definition 3.1.}
Let $S$ be a simple closed curve in the phase plane not passing
through a singular point of the system
    $$
    \dot{x}=P(x,y), \quad \dot{y}=Q(x,y),
    \eqno(3.1)
    $$
where $P(x,y)$ and $Q(x,y)$ are continuous functions (for example,
polynomials), and $M$ be some point on $S.$ If the point $M$ goes
around the curve $S$ in the positive direction (counterclockwise)
one time, then the vector coinciding with the direction of a tangent
to the trajectory passing through the point $M$ is rotated through
the angle $2\pi j$ $(j=0,\pm1,\pm2,\ldots).$ The integer $j$ is
called the \emph{Poincar\'{e} index} of the closed curve $S$
relative to the vector field of system~(3.1) and has the
expression
    $$
    j=\frac{1}{2\pi}\oint_S\frac{P~dQ-Q~dP}{P^2+Q^2}\,.\\[-2mm]
    \eqno(3.2)
    $$
    \par
According to this definition, the index of a node or a focus,
or a center is equal to $+1$ and the index of a saddle is $-1.$
    \par
    \medskip
    \textbf{Theorem 3.1 (First Poincar\'{e} Index Theorem).}
    \emph{If $N,$ $N_f,$ $N_c,$ and $C$ are respectively the number
of nodes, foci, centers, and saddles in a finite part of the phase
plane and $N'$ and $C'$ are the number of nodes and saddles at
infinity, then it is valid the formula}
    $$
    N+N_f+N_c+N'=C+C'+1.
    \eqno(3.3)
    $$
    \par
    \textbf{Theorem 3.2 (Second Poincar\'{e} Index Theorem).}
    \emph{If all singular points are simple, then along an isocline
without multiple points lying in a Poincar\'{e} hemisphere which
is obtained by a stereographic projection of the phase plane,
the singular points are distributed so that a saddle is followed
by a node or a focus, or a center and vice versa. If two points
are separated by the equator of the Poincar\'{e} sphere, then
a saddle will be followed by a saddle again and a node or
a focus, or a center will be followed by a node or a focus,
or a center.}
    \medskip
    \par
Consider system (1.3). Its finite singularities are determined
by the algebraic system
    $$
    x\,(1+\beta_{1}\,x+\ldots+\beta_{2l}\,x^{2l})=0, \quad y=0.
    \eqno(3.4)
    $$
It always has an anti-saddle at the origin and, in general, can have at most
$2l+1$ finite singularities which lie on the $x$-axis and, according to
Theorem~3.2, are distributed so that a saddle (or~saddle-node) is followed
by a node or a focus, or a center and vice versa. At~in\-fi\-nity, system~(1.3)
has two singular points: a node at the ``ends'' of the $x$-axis and a saddle
at the ``ends'' of the $y$-axis. For studying the infinite singularities,
the methods applied in \cite{BL} for Rayleigh's and van der Pol's equations
and also Erugin's two-isocline method developed in \cite{G} can be used
(see \cite{g08a}, \cite{g09a}).
    \par
Following \cite{G}, we will study limit cycle bifurcations of~(1.3) by means of
a~canonical system containing field rotation parameters of~(1.3)~\cite{BL},~\cite{G}.
    \par
    \medskip
\noindent\textbf{Theorem 3.3.}
    \emph{The Li\'{e}nard polynomial system $(1.3)$ with limit cycles can be reduced
to the canonical form}
    \vspace{-4mm}
    $$
    \begin{array}{c}
    \dot{x}=y \equiv P(x,y),\\[1mm]
    \dot{y}=-x\,(1+\beta_{1}\,x\pm\,x^{2}+\ldots+\beta_{2l-1}\,x^{2l-1}\pm\,x^{2l})\\
    +\,y\,(\alpha_{0}+x+\alpha_{2}\,x^{2}+\ldots+x^{2k-1}+\alpha_{2k}\,x^{2k}) \equiv Q(x,y),
    \end{array}
    \eqno(3.5)
    $$
\emph{where $\beta_{1},$ $\beta_{3},\ldots,$ $\beta_{2l-1}$ are fixed and $\alpha_{0},$
$\alpha_{2},\ldots,$ $\alpha_{2k}$ are field rotation para\-me\-ters of $(3.5).$}
    \medskip
    \par
\noindent\textbf{Proof.} Let all the parameters $\alpha_{i},$ $i=0,1,\ldots,2k,$ vanish
in system~(3.5),
    $$
    \dot{x}=y,
    \quad
    \dot{y}=-x\,(1+\beta_{1}\,x+\beta_{2}\,x^{2}+\ldots+\beta_{2l}\,x^{2l}),
    \eqno(3.6)
    $$
and consider the corresponding equation
    $$
    \frac{dy}{dx}
    =\frac{-x\,(1+\beta_{1}\,x+\beta_{2}\,x^{2}+\ldots+\beta_{2l}\,x^{2l})}{y}
    \equiv F(x,y).
    \eqno(3.7)
    $$
Since $F(x,-y)=-F(x,y),$ the direction field of (3.7) (and the vector field of (3.6) as
well) is symmetric with respect to the $x$-axis. It follows that for arbitrary values of
the parameters $\beta_{j},$ $j=1,2,\ldots,2l,$ system~(3.6) has centers as anti-saddles and
cannot have limit cycles surrounding these points. Therefore, without loss of generality,
all the even parameters $\beta_{j}$ of system~(1.3) can be supposed to be equal, e.\,g.,
to $\pm1$: $\beta_{2}=\beta_{4}=\beta_{6}=\ldots=\pm1.$
    \par
Let now all the parameters $\alpha_{i}$ with even indexes and $\beta_{j}$ with odd indexes
vanish in system~(3.5),
    $$
    \dot{x}=y,
    \quad
    \dot{y}=-x\,(1\pm\,x^{2}\pm\ldots\pm\,x^{2l})
    +y\,(\alpha_{1}\,x+\alpha_{3}\,x^{3}+\ldots+\alpha_{2k-1}\,x^{2k-1}),
    \eqno(3.8)
    $$
and consider the corresponding equation
    $$
    \begin{array}{c}
    \displaystyle\frac{dy}{dx}
    =\frac{-x\,(1\pm\,x^{2}\pm\ldots\pm\,x^{2l})
    +y\,(\alpha_{1}\,x+\alpha_{3}\,x^{3}+\ldots+\alpha_{2k-1}\,x^{2k-1})}{y}\\
    \equiv G(x,y).
    \end{array}
    \eqno(3.9)
    $$
Since $G(-x,y)=-G(x,y),$ the direction field of~(3.9) (and the
vector field of (3.8) as well) is symmetric with respect to the
$y$-axis. It follows that for arbitrary values of the parameters
$\alpha_{1},$ $\alpha_{3},\ldots,$ $\alpha_{2k-1}$ system~(3.6)
has centers as anti-saddles and cannot have limit cycles surrounding
these points. Therefore, without loss of generality, all the odd
parameters $\alpha_{i}$ of system~(1.3) can be supposed to be equal,
e.\,g., to~1: $\alpha_{1}=\alpha_{3}=\ldots=\alpha_{2k-1}=1.$
    \par
Inputting the odd parameters $\beta_{1},$ $\beta_{3},\ldots,$ $\beta_{2l-1}$
into system~(3.8),
    $$
    \begin{array}{c}
    \dot{x}=y\equiv R(x,y),\\[1mm]
    \dot{y}=-x\,(1+\beta_{1}\,x\pm\,x^{2}+\beta_{3}\,x^{3}
    \pm\,x^{4}+\ldots+\beta_{2l-1}\,x^{2l-1}\pm\,x^{2l})\\
    +\,y\,(x+x^{3}+\ldots+x^{2k-1})\equiv S(x,y),\\[1mm]
    \end{array}
    \eqno(3.10)
    $$
and calculating the determinants
    $$
\Delta_{\beta_{1}}=RS'_{\beta_{1}}-SR'_{\beta_{1}}=-x^{2}y,
    $$
    $$
\Delta_{\beta_{3}}=RS'_{\beta_{3}}-SR'_{\beta_{3}}=-x^{4}y,
    $$
    $$.\,.\,.\,.\,.\,.\,.\,.\,.\,.\,.\,.\,.\,.\,
    \,.\,.\,.\,.\,.\,.\,.\,.\,.\,.\,.\,.\,.\,.\,.
    $$
    $$
\Delta_{\beta_{2l-1}}=RS'_{\beta_{2l-1}}-SR'_{\beta_{2l-1}}=-x^{2l}y,\\[2mm]
    $$
    \par
\noindent we can see that the vector field of~(3.10) is rotated symmetrically (in opposite directions)
with respect to the $x$-axis and that the finite singularities (centers and saddles) of~(3.10)
moving along the $x$-axis (except the center at the origin) do not change their type or join
in saddle-nodes. Therefore, we can fix the odd parameters $\beta_{1},$ $\beta_{3},\ldots,$
$\beta_{2l-1}$ in system~(3.5), fixing the position of its finite singularities on the $x$-axis.
    \par
To prove that the even parameters $\alpha_{0},$ $\alpha_{2},\ldots,$ $\alpha_{2k}$
rotate the vector field of~(3.5), let us calculate the following determinants:
    $$
\Delta_{\alpha_{0}}=PQ'_{\alpha_{0}}-QP'_{\alpha_{0}}=y^{2}\geq0,
    $$
    $$
\Delta_{\alpha_{2}}=PQ'_{\alpha_{2}}-QP'_{\alpha_{2}}=x^{2}y^{2}\geq0,
    $$
    $$.\,.\,.\,.\,.\,.\,.\,.\,.\,.\,.\,.\,.\,.\,.\,.
    \,.\,.\,.\,.\,.\,.\,.\,.\,.\,.\,.\,.\,.\,.\,.\,.\,.
    $$
    $$
\Delta_{\alpha_{2k}}=PQ'_{\alpha_{2k}}-QP'_{\alpha_{2k}}=x^{2k}y^{2}\geq0.\\[2mm]
    $$
    \par
By definition of a field rotation parameter~\cite{BL},~\cite{G}, for increasing
each of the parameters $\alpha_{0},$ $\alpha_{2},\ldots,$ $\alpha_{2k},$
under the fixed others, the vector field of system~(3.5) is rotated in
the positive direction (counter\-clock\-wise) in the whole phase plane;
and, conversely, for decreasing each of these parameters, the vector
field of~(3.5) is rotated in the negative direction (clock\-wise).
    \par
Thus, for studying limit cycle bifurcations of~(1.3), it is sufficient
to consider the canonical system~(3.5) containing only its even parameters
$\alpha_{0},$ $\alpha_{2},\ldots,$ $\alpha_{2k}$ which rotate the vector
field of~(3.5) under the fixed others. The theorem is proved.\qquad $\Box$
    \par
    \medskip
By means of the canonical system~(3.5), let us study global limit cycle bifurcations of~(1.3)
and prove the following theorem.
    \par
    \medskip
\noindent\textbf{Theorem 3.4.}
    \emph{The general Li\'{e}nard polynomial system $(1.3)$ can have at most $k+l$~limit cycles,
$k$~surrounding the origin and $l$~surrounding one by one the other singularities of $(1.3).$}
    \medskip
    \par
\noindent\textbf{Proof.} According to Theorem~3.3, for the study of limit cycle bifurcations
of system~(1.3), it is sufficient to consider the canonical system~(3.5) containing the field
rotation parameters $\alpha_{0},$ $\alpha_{2},\ldots,$ $\alpha_{2k}$ of~(1.3) under the fixed
its parameters $\beta_{1},$ $\beta_{3},\ldots,$ $\beta_{2l-1}.$
    \par
Let all these parameters vanish:
    $$
    \dot{x}=y,
    \quad
    \dot{y}=-x\,(1\pm\,x^{2}\pm\ldots\pm\,x^{2l})+y\,(x+x^{3}+\ldots+x^{2k-1}).
    \eqno(3.11)
    $$
System (3.11) is symmetric with respect to the $y$-axis and has centers as anti-saddles.
Its center domains are bounded by either separatrix loops or digons of the saddles of~(3.11)
lying on the $x$-axis. If to input the parameters $\beta_{1},$ $\beta_{3},\ldots,$ $\beta_{2l-1}$
into~(3.11) successively, we will get again system~(3.10) the vector field of which is rotated
symmetrically (in opposite directions) with respect to the $x$-axis. The finite singularities
(centers and saddles) of~(3.10) moving along the $x$-axis (except the center at the origin)
do not change their type or join in saddle-nodes and the center domains will be bounded
by separatrix loops of the saddles (or saddle-nodes) of~(3.10)~\cite{BL},~\cite{G}.
    \par
Let us input successively the field rotation parameters $\alpha_{0},$ $\alpha_{2},\ldots,$ $\alpha_{2k}$
into system~(3.10) beginning with the parameters at the highest degrees of~$x$ and alternating with their
signs (see \cite{g08a}, \cite{g09a}). So, begin with the parameter~$\alpha_{2k}$ and let, for definiteness,
$\alpha_{2k}>0\!:$
    \vspace{-1mm}
    $$
    \begin{array}{c}
    \dot{x}=y,\\[1mm]
    \dot{y}=-x\,(1+\beta_{1}\,x\pm\,x^{2}+\beta_{3}\,x^{3}
    \pm\,x^{4}+\ldots+\beta_{2l-1}\,x^{2l-1}\pm\,x^{2l})\\
    +\,y\,(x+x^{3}+\ldots+x^{2k-1}+\alpha_{2k}\,x^{2k}).\\[1mm]
    \end{array}
    \eqno(3.12)
    $$
In this case, the vector field of~(3.12) is rotated in the positive direction (counterclockwise) turning
the center at the origin into a nonrough (weak) unstable focus. All the other centers become rough unstable
foci~\cite{BL},~\cite{G}.
    \par
Fix $\alpha_{2k}$ and input the parameter $\alpha_{2k-2}<0$ into~(3.12):
    \vspace{-1mm}
    $$
    \begin{array}{c}
    \dot{x}=y,\\[1mm]
    \dot{y}=-x\,(1+\beta_{1}\,x\pm\,x^{2}+\beta_{3}\,x^{3}
    \pm\,x^{4}+\ldots+\beta_{2l-1}\,x^{2l-1}\pm\,x^{2l})\\
    +\,y\,(x+x^{3}+\ldots+\alpha_{2k-2}x^{2k-2}+x^{2k-1}+\alpha_{2k}\,x^{2k}).\\[1mm]
    \end{array}
    \eqno(3.13)
    $$
Then the vector field of~(3.13) is rotated in the opposite direction
(clockwise) and the focus at the origin immediately changes the character
of its stability (since its degree of nonroughness decreases and the sign
of the field rotation parameter at the lower degree of~$x$ changes)
generating a stable limit cycle. All the other foci will also generate
stable limit cycles for some values of $\alpha_{2k-2}$ after changing
the character of their stability. Under further decreasing $\alpha_{2k-2},$
all the limit cycles will expand disappearing on separatrix cycles
of~(3.13)~\cite{BL},~\cite{G}.
    \par
Denote the limit cycle surrounding the origin by $\Gamma\!_{1},$ the domain
outside the cycle by $D_{1},$ the domain inside the cycle by $D_{2}$ and
consider logical possibilities of the appearance of other (semi-stable)
limit cycles from a ``trajectory concentration'' surrounding this singular
point. It is clear that, under decreasing the parameter $\alpha_{2k-2},$
a semi-stable limit cycle cannot appear in the domain $D_{2},$ since
the focus spirals filling this domain will untwist and the distance
between their coils will increase because of the vector field
rotation~\cite{g08a},~\cite{g09a}.
    \par
By contradiction, we can also prove that a semi-stable limit cycle
cannot appear in the domain $D_{1}.$ Suppose it appears in this
domain for some values of the parameters $\alpha_{2k}^{*}>0$ and
$\alpha_{2k-2}^{*}<0.$ Return to system (3.10) and change the
inputting order for the field rotation parameters. Input first
the parameter $\alpha_{2k-2}<0\!:$
    \vspace{-1mm}
    $$
    \begin{array}{c}
    \dot{x}=y,\\[1mm]
    \dot{y}=-x\,(1+\beta_{1}\,x\pm\,x^{2}+\beta_{3}\,x^{3}
    \pm\,x^{4}+\ldots+\beta_{2l-1}\,x^{2l-1}\pm\,x^{2l})\\
    +\,y\,(x+x^{3}+\ldots+\alpha_{2k-2}x^{2k-2}+x^{2k-1}).\\[1mm]
    \end{array}
    \eqno(3.14)
    $$
Fix it under $\alpha_{2k-2}=\alpha_{2k-2}^{*}.$ The vector field
of~(3.14) is rotated clockwise and the origin turns into a nonrough
stable focus. Inputting the parameter $\alpha_{2k}>0$ into (3.14),
we get again system (3.13) the vector field of which is rotated
counterclockwise. Under this rotation, a stable limit cycle $\Gamma\!_{1}$
will appear from a separatrix cycle for some value of $\alpha_{2k}.$
This cycle will contract, the outside spirals winding onto the cycle
will untwist and the distance between their coils will increase under
increasing $\alpha_{2k}$ to the value $\alpha_{2k}^{*}.$ It follows that
there are no values of $\alpha_{2k-2}^{*}<0$ and $\alpha_{2k}^{*}>0$
for which a semi-stable limit cycle could appear in the domain~$D_{1}.$
    \par
This contradiction proves the uniqueness of a limit cycle surrounding
the origin in system~(3.13) for any values of the parameters $\alpha_{2k-2}$
and $\alpha_{2k}$ of different signs. Obviously, if these parameters have
the same sign, system~(3.13) has no limit cycles surrounding the origin at~all.
On the same reason, this system cannot have more than $l$~limit cycles surrounding
the other singularities (foci or nodes) of~(3.13) one by one.
    \par
Let system~(3.13) have the unique limit cycle $\Gamma\!_{1}$ surrounding the origin
and $l$~limit cycles surrounding the other antisaddles of~(3.13). Fix the parameters
$\alpha_{2k}>0,$ $\alpha_{2k-2}<0$ and input the third parameter, $\alpha_{2k-4}>0,$
into this system:
    \vspace{-2mm}
    $$
    \begin{array}{c}
    \dot{x}=y,\\[1mm]
    \dot{y}=-x\,(1+\beta_{1}\,x\pm\,x^{2}+\beta_{3}\,x^{3}
    \pm\,x^{4}+\ldots+\beta_{2l-1}\,x^{2l-1}\pm\,x^{2l})\\
    +\,y\,(x+x^{3}+\ldots+\alpha_{2k-4}x^{2k-4}
    +\alpha_{2k-2}x^{2k-2}+x^{2k-1}
    +\alpha_{2k}\,x^{2k}).\\[1mm]
    \end{array}
    \eqno(3.15)
    $$
The vector field of~(3.15) is rotated counterclockwise, the focus
at the origin changes the character of its stability and the second
(unstable) limit cycle, $\Gamma_{2},$ immediately appears from this point.
The limit cycles surrounding the other singularities of~(3.15) can only disappear
in the corresponding foci (because of their roughness) under increasing the parameter
$\alpha_{2k-4}.$ Under further increasing $\alpha_{2k-4},$ the limit cycle $\Gamma_{2}$
will join with $\Gamma\!_{1}$ forming a semi-stable limit cycle, $\Gamma_{\!12},$ which
will disappear in a ``trajectory concentration'' surrounding the origin. Can another
semi-stable limit cycle appear around the origin in addition to $\Gamma_{\!12}?$
It is clear that such a limit cycle cannot appear either in the domain $D_{1}$
bounded on the inside by the cycle $\Gamma\!_{1}$ or in the domain $D_{3}$
bounded by the origin and $\Gamma_{2}$ because of the increasing distance
between the spiral coils filling these domains under increasing the
parameter~$\alpha_{2k-4}$~\cite{g08a},~\cite{g09a}.
    \par
To prove the impossibility of the appearance of a semi-stable limit
cycle in the domain $D_{2}$ bounded by the cycles $\Gamma\!_{1}$
and $\Gamma_{2}$ (before their joining), suppose the contrary, i.\,e.,
that for some set of values of the parameters, $\alpha_{2k}^{*}>0,$
$\alpha_{2k-2}^{*}<0,$ and $\alpha_{2k-4}^{*}>0,$ such a semi-stable
cycle exists. Return to system~(3.10) again and input first the parameters
$\alpha_{2k-4}>0$ and $\alpha_{2k}>0\!:$
    \vspace{-2mm}
    $$
    \begin{array}{c}
    \dot{x}=y,\\[1mm]
    \dot{y}=-x\,(1+\beta_{1}\,x\pm\,x^{2}+\beta_{3}\,x^{3}
    \pm\,x^{4}+\ldots+\beta_{2l-1}\,x^{2l-1}\pm\,x^{2l})\\
    +\,y\,(x+x^{3}+\ldots+\alpha_{2k-4}x^{2k-4}+x^{2k-3}
    +\alpha_{2k}\,x^{2k}).\\[1mm]
    \end{array}
    \eqno(3.16)
    $$
Both parameters act in a similar way: they rotate the vector field of~(3.16)
counterclockwise turning the origin into a nonrough unstable focus.
    \par
Fix these parameters under $\alpha_{2k-4}=\alpha_{2k-4}^{*},$
$\alpha_{2k}=\alpha_{2k}^{*}$ and input the parameter $\alpha_{2k-2}<0$
into~(3.16) getting again system~(3.15). Since, by our assumption, this system
has two limit cycles surrounding the origin for $\alpha_{2k-2}>\alpha_{2k-2}^{*},$
there exists some value of the parameter, $\alpha_{2k-2}^{12}$
$(\alpha_{2k-2}^{*}<\alpha_{2k-2}^{12}<0),$ for which a semi-stable
limit cycle, $\Gamma\!_{12},$ appears in system~(3.15) and then
splits into a stable cycle, $\Gamma\!_{1},$ and an unstable cycle,
$\Gamma_{2},$ under further decreasing $\alpha_{2k-2}.$ The formed
domain $D_{2}$ bounded by the limit cycles $\Gamma\!_{1},$ $\Gamma_{2}$
and filled by the spirals will enlarge since, on the properties of
a field rotation parameter, the interior unstable limit cycle $\Gamma_{2}$
will contract and the exterior stable limit cycle $\Gamma\!_{1}$ will expand
under decreasing $\alpha_{2k-2}.$ The distance between the spirals of the domain
$D_{2}$ will naturally increase, which will prevent the appearance of a semi-stable
limit cycle in this domain for $\alpha_{2k-2}<\alpha_{2k-2}^{12}$~\cite{g08a},~\cite{g09a}.
    \par
Thus, there are no such values of the parameters, $\alpha_{2k}^{*}>0,$
$\alpha_{2k-2}^{*}<0,$ and $\alpha_{2k-4}^{*}>0,$ for which system~(3.15)
would have an additional semi-stable limit cycle surrounding the origin.
Obviously, there are no other values of the para\-meters $\alpha_{2k},$
$\alpha_{2k-2},$ and $\alpha_{2k-4}$ for which system~(3.15) would have
more than two limit cycles surrounding this singular point. On the same reason,
additional semi-stable limit cycles cannot appear around the other singularities
(foci or nodes) of~(3.15). Therefore, $2+l$ is the maximum number of limit cycles
in system~(3.15).
    \par
Suppose that system (3.15) has two limit cycles, $\Gamma\!_{1}$ and $\Gamma_{2},$
surrounding the origin and $l$~limit cycles surrounding the other antisaddles of~(3.15)
(this is always possible if $\alpha_{2k}\gg-\alpha_{2k-2}\gg\alpha_{2k-4}>0).$
Fix the parameters $\alpha_{2k},$ $\alpha_{2k-2},$ $\alpha_{2k-4}$ and consider
a more general system inputting the fourth parameter, $\alpha_{2k-6}<0,$
into~(3.15):
    \vspace{-2mm}
    $$
    \begin{array}{c}
    \dot{x}=y,\\[1mm]
    \dot{y}=-x\,(1+\beta_{1}\,x\pm\,x^{2}+\beta_{3}\,x^{3}
    \pm\,x^{4}+\ldots+\beta_{2l-1}\,x^{2l-1}\pm\,x^{2l})\\
    +\,y\,(x+x^{3}+\ldots+\alpha_{2k-6}x^{2k-6}+x^{2k-5}
    +\ldots+\alpha_{2k}\,x^{2k}).\\[1mm]
    \end{array}
    \eqno(3.17)
    $$
For decreasing $\alpha_{2k-6},$ the vector field of~(3.17) will
be rotated clockwise and the focus at the origin will immediately
change the character of its stability generating a third
(stable) limit cycle, $\Gamma_{3}.$ With further decreasing
$\alpha_{2k-6},$ $\Gamma_{3}$ will join with $\Gamma_{2}$ forming
a semi-stable limit cycle, $\Gamma_{\!23},$ which will disappear
in a ``trajectory concentration'' surrounding the origin; the cycle
$\Gamma\!_{1}$ will expand disappearing on a separatrix cycle of~(3.17).
    \par
Let system (3.17) have three limit cycles surrounding the origin:
$\Gamma\!_{1},$ $\Gamma_{2},$ $\Gamma_{3}.$ Could an additional
semi-stable limit cycle appear with decreasing $\alpha_{2k-6}$
after splitting of which system~(3.17) would have five limit cycles
around the origin? It is clear that such a limit cycle cannot appear
either in the domain $D_{2}$ bounded by the cycles $\Gamma\!_{1}$ and
$\Gamma_{2}$ or in the domain $D_{4}$ bounded by the origin and
$\Gamma_{3}$ because of the increasing distance between the spiral
coils filling these domains after decreasing $\alpha_{2k-6}.$
Consider two other domains: $D_{1}$ bounded on the inside by the
cycle $\Gamma\!_{1}$ and $D_{3}$ bounded by the cycles
$\Gamma_{2}$ and $\Gamma_{3}.$ As before, we will prove
the impossibility of the appearance of a semi-stable limit
cycle in these domains by contradiction.
    \par
Suppose that for some set of values of the parameters
$\alpha_{2k}^{*}>0,$ $\alpha_{2k-2}^{*}<0,$ $\alpha_{2k-4}^{*}>0,$
and $\alpha_{2k-6}^{*}<0$ such a semi-stable cycle exists. Return
to system~(3.10) again, input first the parameters $\alpha_{2k-6}<0,$
$\alpha_{2k-2}<0$ and then the parameter $\alpha_{2k}>0\!:$
    \vspace{-2mm}
    $$
    \begin{array}{c}
    \dot{x}=y,\\[1mm]
    \dot{y}=-x\,(1+\beta_{1}\,x\pm\,x^{2}+\beta_{3}\,x^{3}
    \pm\,x^{4}+\ldots+\beta_{2l-1}\,x^{2l-1}\pm\,x^{2l})\\
    +\,y\,(x+x^{3}+\ldots+\alpha_{2k-6}x^{2k-6}+\ldots
    +\alpha_{2k-2}x^{2k-2}+x^{2k-3}+\alpha_{2k}\,x^{2k}).\\[1mm]
    \end{array}
    \eqno(3.18)
    $$
Fix the parameters $\alpha_{2k-6},$ $\alpha_{2k-2}$ under the values
$\alpha_{2k-6}^{*},$ $\alpha_{2k-2}^{*},$ respectively. With increasing
$\alpha_{2k},$ a separatrix cycle formed around the origin will generate
a stable limit cycle, $\Gamma\!_{1}.$ Fix $\alpha_{2k}$ under the value
$\alpha_{2k}^{*}$ and input the parameter $\alpha_{2k-4}>0$ into~(3.18)
getting system~(3.17).
    \par
Since, by our assumption, (3.17) has three limit cycles for
$\alpha_{2k-4}<\alpha_{2k-4}^{*},$ there exists some value of the
parameter $\alpha_{2k-4}^{23}$ $(0<\alpha_{2k-4}^{23}<\alpha_{2k-4}^{*})$
for which a semi-stable limit cycle, $\Gamma_{23},$ appears in
this system and then splits into an unstable cycle, $\Gamma_{2},$
and a stable cycle, $\Gamma_{3},$ with further increasing
$\alpha_{2k-4}.$ The formed domain $D_{3}$ bounded by the limit
cycles $\Gamma_{2},$ $\Gamma_{3}$ and also the domain $D_{1}$
bounded on the inside by the limit cycle $\Gamma\!_{1}$ will
enlarge and the spirals filling these domains will untwist
excluding a possibility of the appearance of a semi-stable
limit cycle there~\cite{g08a},~\cite{g09a}.
    \par
All other combinations of the parameters $\alpha_{2k},$ $\alpha_{2k-2},$
$\alpha_{2k-4},$ and $\alpha_{2k-6}$ are considered in a similar way.
It follows that system~(3.17) can have at most $3+l$ limit cycles.
    \par
If we continue the procedure of successive inputting the even parameters,
$\alpha_{2k},\ldots,$ $\alpha_{2},$ $\alpha_{0},$ into system (3.10),
it is possible first to obtain $k$~limit cycles surrounding the origin
$(\alpha_{2k}\gg-\alpha_{2k-2}\gg\alpha_{2k-4}\gg-\alpha_{2k-6}\gg\alpha_{2k-8}\gg\ldots)$
and then to conclude that the canonical system (3.5) (i.\,e., the Li\'{e}nard polynomial
system (1.3) as well) can have at most $k+l$~limit cycles, $k$~surrounding the origin
and $l$~surrounding one by one the antisaddles (foci or nodes) of~(3.5) (and~(1.3)
as~well). The theorem is proved.\qquad $\Box$

\end{document}